\documentclass{article}
\usepackage[T2A]{fontenc}
\usepackage[english]{babel}
\usepackage[tbtags]{amsmath}
\usepackage[dvips]{graphicx}
\usepackage{amsfonts,amssymb,mathrsfs,amscd,amsmath,comment}

\date{}

\begin{document}

\author{B.\,S.~Safin\footnote{Nizhnij Novgorod, Russia; E-mail: {\tt boris-nn@yandex.ru}}}
\title {Studying solutions to Diophantine equations using the table of $N$th digital roots of integers}

\maketitle
\begin{abstract}
In this note we recall the definition of the digital root, and apply the notion of the digital root to searching solutions of
Diophantine equations. A~table of arithmetic operations with digital roots is given. This method is
incapable of obtaining complete solutions of equations, but it is frequently useful in determining
whether this equation is solvable in integers or not. A~minor extension of Fermat's little theorem is put forward.
\end{abstract}

\section{Introduction} The concept of the digital root was introduced by Gardner in the book~\cite{Gar}. He defines the digital
root as follows: if all the digits in a given number are added, then the digits
in the sum added, and this continued until only a single digit
remains, that digit is known as the digital root of the original
number. Thus, in the decimal number system, we may obtain the digits 1, 2, 3, 4, 5, 6, 7, 8,~9.
For the time being, we shall be concerned only with positive numbers. For example, the digital root of  123\,456\,789 is calculated as follows:
$1+2+3+4+5+6+7+8+9$, $45=4+5=9$. The digital root of 888 is as follows: $8+8+8=24=2+4=6$, $100=1+0+0=1$, and so~on.
We split all the integer numbers into 9~classes: $1+9k$, $2+9k$, $3+9k$, $4+9k$, $5+9k$, $6+9k$, $7+9k$, $8+9k$, $9+9k$.
The corresponding digital roots of each of the successive classes are 1,      2,       3,       4,      5,      6,      7,       8,~9 for any~$k$.

Now let us compose an infinite table, in which the digital roots of any of the classes will be raised to the $n$th power (see Table~1).
The left column contains numbers $X^n$ of the form $1+9k$, $2+9k$, $3+9k$, $4+9k$, $5+9k$, $6+9k$, $7+9k$, $8+9k$, $9+9k$ raised to the $n$th power.
The rows contain the corresponding digital roots of these numbers.

\begin{figure}[ht]
\begin{center}
\centerline{Table 1}
\includegraphics[width=0.75\textwidth]{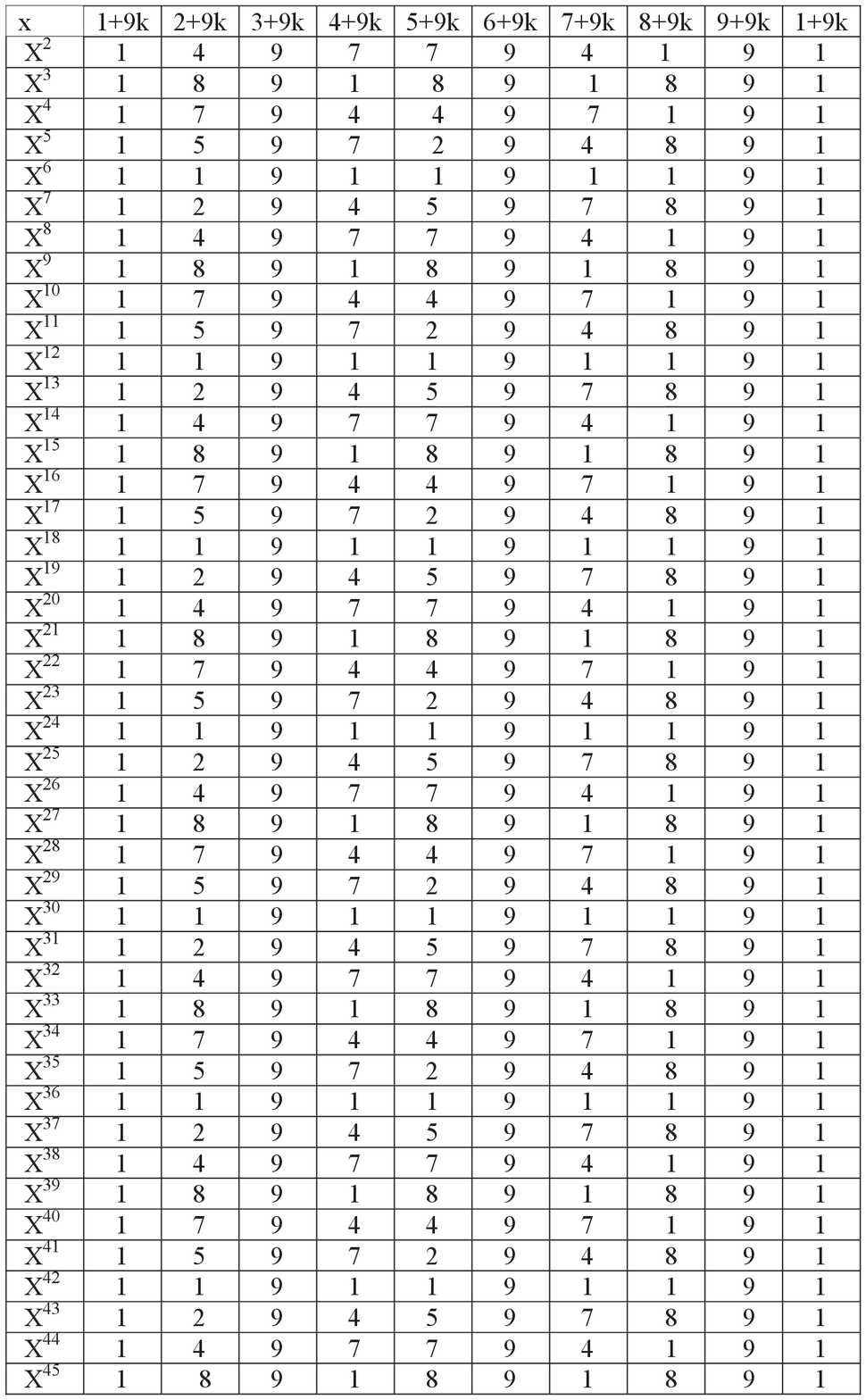}
\end{center}
\end{figure}

Table 2 illustrates the arithmetics of the digital roots (except for the division). Consider, for example, $X^n= (1+9k)^n$.
Assume that $n=2$. Then the digital root of the expression
$(1+9k)^2=1+18k+81k^2$ equals~$1$. We have $ 1+(1+8)k+(8+1)k^2=1+9k+9k^2$, since multiplying any number by the digital root of~$9$
produces the digital root of~$9$ (see Table~2), and hence the digital root of $1+9k+9k^2$ is $1+9+9=19=1+9=1$.
From the viewpoint of the theory  of congruences, the digital root of any number from the 9~classes into which we split all the integer numbers
is a~congruence modulo~9; that is, this is the remainder after dividing by~9.
Therefore, there is no need to raise to a~power all numbers of the form $(m+9k)$, where $m$~ranges from 1 to~9, but instead one
needs to raise to a~power only the remainder after dividing by~9, that is, the digital root of~$m$.

\section {Study of Diophantine equations}

\textbf{1.} The Pell equation $x^2-dy^2=1$, where $d$ is a~non-square natural number.
We rewrite this equation as follows: $x^2=1 + dy^2$, and consider the second row of Table~1, which contains
all digital roots of all square numbers. It is seen from this table that the digital roots of $x^2$, $y^2$ may assume only
four values: 1, 4, 7,~9.
Consider the first case, when the digital root of a~number $x^2$ assumes the value~$1$. In this case, the Pelle equation
reads as follows $1=1+ dy^2$. In order to keep the equality, it is required that the digital root of the expression
$dy^2$ be~$9$. We thus have $1=1+9=10=1+0=1$. In the expression $dy^2$ the digital root  of $y^2$ may assume the values 1, 4, 7,~9.
If the digital root of~$y^2$ is~$1$, then the number~$d$ (and its digital root) may take on the value~$9$; that is,
$d$~may only be of the form~$9k$, since the equality $1=1+9\times 1=1+9=10=
1+0=1$ is preserved only in this case. If $y^2$ assumes other values, for example, 4,~7, then the number~$d$ itself (and its digital root)
may take on only the value~$9k$, since by multiplying the digital roots 4,~7 by the digital root~9 we obtain the digital root~$9$
(see Table~2), the equality $1=1+9=1$ preserving. If the digital root of~$y^2$ assumes the value~9, then $d$ may be arbitrary,
because multiplying by~9 always gives the digital root~$9$ and the equality $1=1+9=1$ is preserved.

If the digital root of $x^2$ assumes other values (4, 7, 9), then the same conclusions can be made about~$d$ along the same lines

\smallskip

\textbf{2. The equation} $ x^2=y^3-2$.  Consider the second and third rows of Table~1. The digital root of $x^2$ may take on the
values 1, 4, 7,~9. The digital root of~$y^3$ may assume the values 1, 8,~9. Substituting the digital root for~$y^3$, we
obtain three variants of possible values for the digital root of~$x^2$: a)~$x^2=9-2=7$, b)~$x^2=8-2=6$, c)~$x^2=1-2=-1$.
Clearly, neither 6 nor~$-1$ ~cannot be the digital root of~$x^2$, inasmuch as $-1=8-9$ or $9-1=8$ (see Table~2), the digital root of~$-1$ is the
digital root of~$8$. So, in order that the negative digital root have positive values one needs to add~9 to~it. In case both roots are
equal, we would have $9-9=9$. Hence, the digital root of $x^2$ may assume the values 6, 7,~8. But it is only~7 that may be the
digital root of a~squared number. Hence, $е$ (see Table~1) may be only of the form ($4+9k)$ or $(5+9k)$, where $k = 0, 1, 2, \dotsc$, and $y$
may only be of the form $(3+9k)$, $(6+9k)$ or~$9k$.

\smallskip

\textbf{3. The Beal conjecture} $A^x+B^y=C^z$. Assume that $A,B,C ,x,y,z$ are natural numbers with $x,y,z>2$. Is it true that
$A,B,C$ have a~common prime factor?
From analysis of Table~1 it is clear that the digital roots of  $A^x$, $B^y$, $C^z$ may assume the following values: 1, 2, 4, 5, 7, 8,~9 for
any $x,y,z>2$. Since the analysis of this equation is cumbersome due to tedious combinatorics, we confine to a~few remarks.
If the power of a~number exceeds~1, then then table shows that there are no numbers whose digital roots are 3 and~6. Hence,
if the digital root of a~number $A^x$ is~1, then the digital root of  $B^y$ may not equal~2, because $1+2=3$. If follows that if
in~$A^x$ one has $A=(1+9k)$, then $x$ is any natural number; if $A=(2+9k)$, then $x=6n$; if $ A=(4+9k)$, then  $x=3n$; if $A=(5+9k)$, then
$x=6n$l if
$A=(7+9k)$,  then $x =3n$, and if $A=(8+9k)$, then  $x=2n$. Hence, since the digital root of $B^y$  may not equal~2,
then $B^y$ may not be of the form $B=(2+9k)$, where $y=(7+6n)$, $B=(5+9k)$, and so $y=5+6n$.
Alos, if the digital root of $A^x$ is~1, then the digital root of~$B^y $ may not equal~$5$, since $1+5=6$ (no such root exists).
Hence, $B^y$ may not be represented as $B=(2+9k)$, where $y=5+6n$, $B=(5+9k)$, where $y=7+6n$.
A~similar analysis applies when the digital root of~$A^x$ is~2, and so~on.

\smallskip

\textbf{4. Fermat's little theorem.}
This theorem reads as follows: if $a$~is a~prime and is not a~multiple of~$p$, where $p$~is a~prime, then
$a^p-a=pn$, where $n$~is an integer.

An analysis of Table~1 yields a small generalization of Fermat's little theorem: if $p$ and $ q$ are odd primes, then $a^p-a^q=3n$, where $n$ is integer.

\begin{figure}[ht]
\begin{center}
\centerline{Table 2. Arithmetics of the digital root}
\smallskip
\centerline{Addition}
\smallskip
\includegraphics[width=0.75\textwidth]{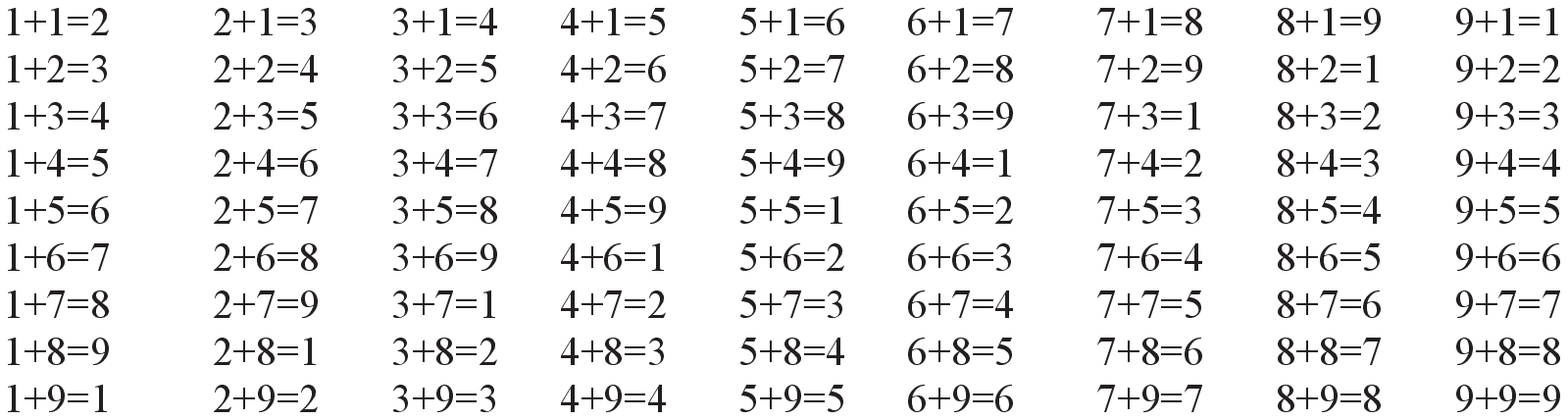}
\smallskip
\centerline{Subtraction}
\smallskip
\includegraphics[width=0.75\textwidth]{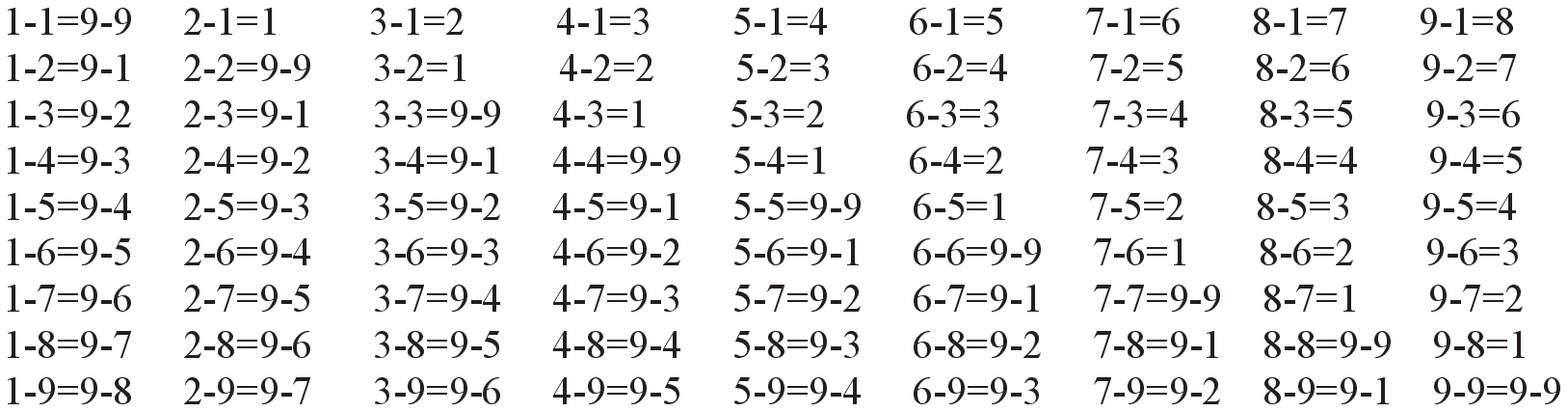}
\smallskip
\centerline{Multiplication}
\smallskip
\includegraphics[width=0.75\textwidth]{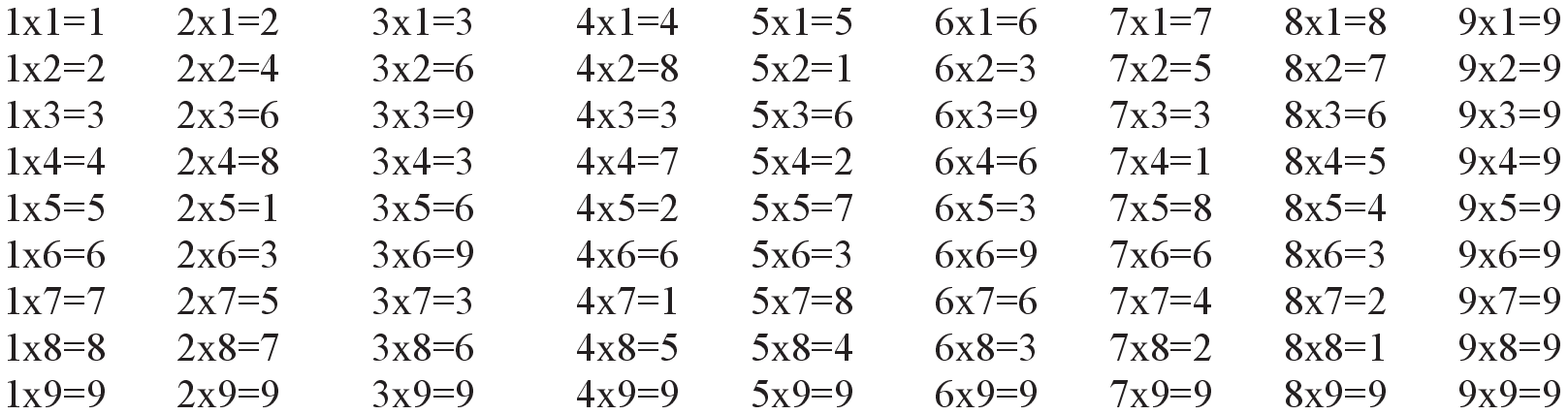}
\end{center}
\end{figure}
\vfill \null

\end{document}